\documentclass[a4paper,preprint,amsmath,pre]{revtex4-1}
\usepackage{graphicx}

\begin{document}
\title{Pattern Formation on Networks with Reactions: A Continuous Time Random Walk Approach}
\author{C. N. Angstmann}
\email{c.angstmann@unsw.edu.au}
\author{I. C. Donnelly}
\email{i.donnelly@unsw.edu.au}
\author{B. I. Henry}
\email{b.henry@unsw.edu.au}
\affiliation{School of Mathematics and Statistics, University of New South Wales, Sydney, NSW 2052, Australia}
\date{\today}

\begin{abstract} 
We derive the generalized master equation for reaction-diffusion on networks from an underlying stochastic process, the continuous time random walk (CTRW). The non-trivial incorporation of the reaction process into the CTRW is achieved by splitting the derivation into two stages. The reactions are treated as birth-death processes and the first stage of the derivation is at the single particle level, taking into account the death process, whilst the second stage considers an ensemble of these particles including the birth process. Using this model we have investigated different types of pattern formation across the vertices on a range of networks. Importantly, the CTRW defines the Laplacian operator on the network in a non \emph{ad-hoc} manner and the pattern formation depends on the structure of this Laplacian. Here we focus attention on CTRWs with exponential waiting times for two cases; one in which the rate parameter is constant for all vertices and the other where the rate parameter is proportional to the vertex degree. This results in nonsymmetric and symmetric CTRW Laplacians respectively. In the case of symmetric Laplacians, pattern formation follows from the Turing instability. However in nonsymmetric Laplacians, pattern formation may be possible with or without a Turing instability. 
\end{abstract}

\maketitle

\section{Introduction}

Networks have been extensively studied as models for highly connected systems in biology \cite{BO2004}, physics \cite{Zallen1983} and the social sciences \cite{Barabasi1999}. Over the past decade there has been a great deal of interest in understanding theoretical properties of transport on networks \cite{LBHS2005,Simonsen2005,TRT2007} with a growing interest in the problem of transport on networks with reactions \cite{GA2006,Nakao2010,Wolfrum2012}.  In this article we provide a detailed derivation of the generalized master equation for transport on networks with reactions. The network model we consider allows for reactions of particles on vertices and diffusion of particles between vertices. A discrete space description of diffusion can be modelled by a random walk  \cite{Bachelier1900, Pearson1905, Einstein1905}. Random walks on networks have been extensively studied in this context \cite{HB1987, NR2004, CBTVK2007}.
Our derivation of the generalized master equation is based on the continuous time random walk (CTRW) formalism in which a random walker waits a time (drawn from a waiting time probability density) before jumping \cite{MW1965, SL1973}. This model has been particularly useful to model diffusion in systems with disorder in waiting times resulting in anomalous diffusion \cite{KS2011}. The CTRW on a spatial continuum or a uniform lattice has been further generalized to include linear reactions  \cite{Sokolov2006, HLW2006}, linear reactions with multiple species \cite{LHW2008} and nonlinear reactions \cite{VR2002, EFJS2008, MFH2010, Fedotov2010, AYL2010}. 

The precise incorporation of reactions into the CTRW model for general networks is nontrivial, as, even in the spatial continuum case, reaction and diffusion processes become entwined \cite{HLW2006, Fedotov2010}. To include reactions in the generalized master equations for CTRWs on a network, we separate the derivation into two stages. The first stage is at the single particle level where the loss of particles due to reactions are treated as a death process. The second stage considers an ensemble of such particles and incorporates the remainder of the reaction-kinetics. The resultant generalized master equation, obtained from an underlying stochastic process, provides a fundamental description of reactions with diffusion on networks. This allows, among other things the study of pattern formation on networks with diffusion and reactions.
For example this model description could be applied to the analysis of diffusion tensor imaging data  \cite{LE2001, Raj2012} and to spatio-temporal models in epidemiology  \cite{Colizza2006, KA2001}. 

In spatial continuum systems reaction-diffusion models form the basis for studies of pattern formation in a wide range of applications. The classic model for pattern formation in these systems is Turing pattern formation \cite{Turing1952} which arises from an instability in the reaction dynamics caused by differing rates of diffusion. Such patterns emerge in biological morphogenesis \cite{GM1972, KM2010, MGM2010, Meinhardt2012, MWBGL2012} chemical reactions \cite{OS1991}, propagation of viruses \cite{SJLL2007, SAMH2012} and ecosystems encompassing competing animals \cite{Mimura1978}. 
There have also been studies of Turing patterns on networks. Significantly, the structure of the network has a direct effect on the resulting pattern \cite{Othmer1971,CPSV2007}. Other studies of pattern formation in networks have considered scale-free networks \cite{GA2004}, coupled reactors \cite{HLM2004}, functional gene networks \cite{DD2006}, multiple coexisting stationary states \cite{Nakao2010, Hata2012, Wolfrum2012, KKM2012}, the effects of feedback \cite{Hata2012} and the formation of traveling fronts \cite{KKM2012}.

In this article we have used the CTRW framework to derive a family of diffusive network Laplacian operators incorporating reactions and we have studied pattern formation in these systems. The reaction-diffusion behaviour with these network Laplacian operators may differ from that with the continuum Laplacian operator \cite{Wardetzky2007}. The generalized master equation that we derived has few restrictions on the form of the waiting time density. Importantly it allows us to model both standard transport, arising from exponential waiting time densities, and anomalous transport, arising from power law waiting time densities. In spatial continuum systems anomalous transport has been shown to alter the onset and nature of Turing patterns \cite{HLW2005,NN2007,HVB2009}. 
However to simplify presentation of results from this model we have confined our further analysis to pattern formation arising from exponential waiting time densities and we have considered both symmetric and nonsymmetric Laplacians in this context. 
We have carried out algebraic analysis and model simulations that show pattern formation on Bar\'abasi-Albert networks and \cite{Barabasi1999}, Watts-Strogatz networks \cite{Watts1998}. Gierer-Meinhardt reaction kinetics were used in these examples as representative of reactions that permit Turing pattern formation on a spatial continuum \cite{GM1972}.
The examples that we considered demonstrate the influence of network topology and network diffusion on pattern formation. 

The remainder of the paper is as follows. In Sec. \ref{sec:derivation} we derive the master equations that describe the CTRW network reaction diffusion model. In Sec. \ref{sec:patterns} we describe different pattern formation mechanisms that may arise depending on the form of the CTRW Laplacian.
In Sec. \ref{sec:numerical} we present numerical simulations of pattern formation in the network models. We conclude with a summary and discussion in Sec. \ref{sec:conclusion}.

\section{Derivation of a CTRW network Reaction Diffusion Master Equations}
\label{sec:derivation}

Diffusion, or diffusive-like phenomena, can arise on a network from a variety of sources depending on the details of the network. In considering diffusion we should begin by defining a stochastic process that makes physical sense for the phenomena on the network. 

A CTRWis a stochastic process that naturally limits to diffusion in the continuum \cite{MW1965, SL1973}. To model a reaction-diffusion process we assume that the motion of each individual particle can be expressed as a CTRW. That is to say the particles will jump from vertex to vertex on the network according to the edges present. On each vertex they will wait for a random time before randomly jumping to a connected vertex. The model of random walks with reactions and variation in node degree is analogous to the reaction-diffusion model with spatially varying diffusion  \cite{PMM2005}.
 
The reactions occur between particles that occupy the same vertex. We consider the reactions to be a birth-death process of the particles. Particles will be created according to some probability, and destroyed according to a different probability. These probabilities may depend on the density of other particles on the vertex. We can then derive the equations that govern the evolution of a single particle in time. In this manner the evolution of a single particle is subject to only the death probabilities. The birth process is included by summing the initial conditions of each single particle CTRW with the probability that a particle was created on a particular vertex at an instant in time. 

The assumption that each vertex is a well mixed system is made so that the rate of the probability of a particle being destroyed in reactions is not dependent on the amount of time a particle has been waiting on the vertex. 
We assume that the number of particles at any given vertex is sufficiently large to justify the well mixed approximation i.e. law of mass action and reaction kinetics.
We also assume that the waiting time for each newly created particle is independent of the waiting times of the parent particles. This is similar to Model B in \cite{Fedotov2010} which was first considered by Vlad and Ross \cite{VR2002}. The first step in our derivation is to obtain the master equation for the evolution of a single particle subject to a probability of death that is inhomogeneous in time and space.

\subsection{Single Particle CTRW Death Process Density}
Consider a network whose vertices form the set $\mathbf{W} = \{w_1,  ...,  w_J \}$ where $J$ is the number of vertices. Let $\rho(w_j,t|w_0,0)$ be the probability density for a random walker to be on vertex $w_j$ at time $t$ given it started on vertex $w_0 \in \mathbf{W}$ at time $t=0$. 

Define $q_n(w_j,t|w_0,0)$ as the conditional probability density for arriving at vertex $w_j$ at time $t$ after $n$ steps. We define the reaction survival function, $e^{-\int\limits_{t'}^t \beta(w_j,t'') dt''}$, as the probability that a particle stays alive from $t'$ to $t$ given it does not leave vertex $w_j$ and $\beta$ is a death rate that in general may depend on vertex and time.
The initial condition for $n=0$ is given by 
\begin{equation}
 q_0(w_j,t|w_0,0) = \delta_{w_j,w_0} \delta(t-0^+).
\label{eqn:bound_cond}
\end{equation}
In general, we can write
\begin{equation}
 q_{n+1}(w_j,t|w_0,0) = \sum\limits_{i=1}^J \int\limits_0^t \Psi(w_i,w_j,t,t') e^{-\int\limits_{t'}^t \beta(w_i,t'') dt''}q_n(w_i,t'|w_0,0) dt'
\label{eqn:q_n+1}
\end{equation}
where $\Psi(w_i,w_j,t,t')$ is the probability density of the transition to vertex $w_j$ at time $t$ given the random walker arrived at vertex $w_i$ at an earlier time $t'$ after $n$ steps.

We assume that $\Psi(w_i,w_j,t,t')$ may be expressed as a product of two independent densities; a jump density $\lambda(w_j,w_i)$ and a waiting time density $\psi(w_i,t-t')$ so that
\begin{equation}
\label{eqn:prob_trans}
 \Psi(w_i,w_j,t,t') = \lambda(w_j,w_i) \psi(w_i,t-t')
\end{equation}
where $\lambda$ and $\psi$ must satisfy the normalizations
\begin{equation}
\sum\limits_{j=1}^J \lambda(w_j,w_i)  =1 \qquad \textrm{ for fixed}\  w_i
\end{equation}
and 
\begin{equation}
\int\limits_{t'}^\infty \psi(w_i,t-t') \, dt =1 \qquad \textrm{ for fixed}\  w_i \ \textrm{and} \ t'.
\end{equation}
The separation in Eq. (\ref{eqn:prob_trans}) facilitates the derivation of the generalized master equations in this paper. The inclusion of the vertex (spatial) dependence in the waiting time is more general than the standard independence assumption used in CTRW derivations \cite{MK2000}. 

The conditional density for the walker to arrive at $w_j$ at time $t$ after any number of steps is found by summing over all $n$ steps using Eqs. (\ref{eqn:q_n+1}) and (\ref{eqn:prob_trans}):
\begin{align}
 q(w_j,t|w_0,0) & = \sum\limits_{n=0}^{\infty} q_n(w_j,t|w_0,0) \nonumber \\
		& = \delta_{w_j,w_0} \delta(t-0^+) + \sum\limits_{n=0}^{\infty} \sum\limits_{i=1}^J \int\limits_0^t \Psi(w_i,w_j,t,t')e^{-\int\limits_{t'}^t \beta(w_i,t'') dt''}q_n(w_i,t'|w_0,0) dt' \nonumber \\
		& = \delta_{w_j,w_0} \delta(t-0^+) + \sum\limits_{i=1}^J \lambda(w_j,w_i) \int\limits_0^t \psi(w_i,t-t') e^{-\int\limits_{t'}^t \beta(w_i,t'') dt''}q(w_i,t'|w_0,0) dt'.
		\label{eqn:gen_q}
\end{align}

We can then define the conditional probability density for the random walker to be at vertex $w_j$ at time $t$;
\begin{equation}
 \rho(w_j,t|w_0,0) =  \int\limits_0^t \phi(w_j,t-t') e^{-\int\limits_{t'}^t \beta(w_j,t'') dt''}q(w_j,t'|w_0,0) dt',
\label{eqn:rho_defn}
\end{equation}
where $\phi(w_j,t-t')$ is the probability that the particle does not jump during the period of time $t-t'$:
\begin{equation}
 \phi(w_j,t-t') = 1 - \int\limits_0^{t-t'} \psi(w_j,t'')dt''.
\label{eqn:surv_prob}
\end{equation}

\subsection{Single Particle CTRW Death Master Equation}
 
The derivation of the master equations describing a CTRW death process on a network is similar to the derivations presented in \cite{CGS2005}, \cite{Angstmann2011} and \cite{Fedotov2010} . Formally, the integrals over probability densities should be treated as Riemann-Stieltjes integrals and care has to be taken due to the discontinuity in the arrival density $q(w_j,t|w_0,0)$ at time $t = 0$ \cite{Angstmann2011}. To do this, write
\begin{equation}
 q(w_j,t|w_0,0) = \delta_{w_j,w_0} \delta(t-0^+) + q^+(w_j,t|w_0,0)
\label{eqn:dis_bc}
\end{equation}
where $q^+$ is right side continuous at $t = 0$.
Thus by substitution of Eq. (\ref{eqn:dis_bc}) into Eq. (\ref{eqn:rho_defn}) we get 
\begin{equation}
 \rho(w_j,t|w_0,0) = \delta_{w_j,w_0} \phi(w_j,t) e^{-\int\limits_{0}^t \beta(w_j,t') dt'} + \int\limits_0^t q^+(w_j,t'|w_0,0) e^{-\int\limits_{t'}^t \beta(w_j,t'') dt''}\phi(w_j,t-t')dt' .
\label{eqn:rho_master_initial}
\end{equation}
We now differentiate this equation with respect to time, using the Leibniz rule for differentiating under the integral sign, to obtain
\begin{eqnarray}
 \frac{d \rho(w_j,t|w_0,0)}{ d t} &=&  q^+(w_j,t|w_0,0)  - \int\limits_0^t q^+(w_j,t'|w_0,0) e^{-\int\limits_{t'}^t \beta(w_j,t'') dt''}\psi(w_j,t-t')dt' \nonumber \\
  & &-  \beta(w_j,t) \rho(w_j,t|w_0,0) - \delta_{w_j,w_0} e^{-\int\limits_{0}^t \beta(w_j,t'') dt''}\psi(w_j,t).
\label{eqn:diff_rho}
\end{eqnarray}
Define the flux leaving vertex $w_j$ at time $t$ as 
\begin{align}
 i(w_j,t|w_0,0) &=  \delta_{w_j,w_0} e^{-\int\limits_0^t \beta(w_j,t')dt'} \psi(w_j,t) + \int\limits_0^t q^+(w_j,t'|w_0,0) e^{-\int\limits_{t'}^t \beta(w_j,t'') dt''}\psi(w_j,t-t')dt' \\
 		& =   \int\limits_0^t q(w_j,t'|w_0,0) e^{-\int\limits_{t'}^t \beta(w_j,t'') dt''}\psi(w_j,t-t')dt'.
\label{eqn:flux_def}
\end{align}

We can then rewrite Eq. (\ref{eqn:diff_rho}) as
\begin{equation}
 \frac{d\rho(w_j,t|w_0,0)}{d t} = q^+(w_j,t|w_0,0) - i(w_j,t|w_0,0) - \beta(w_j,t)\rho(w_j,t|w_0,0).
\label{eqn:diff_rho_flux}
\end{equation}

Using Eqs. (\ref{eqn:gen_q}), (\ref{eqn:dis_bc}) and (\ref{eqn:flux_def}), the rate of arrivals at vertex $w_j$ can be expressed as 
\begin{equation}
 q^+(w_j,t|w_0,0) = \sum_{i=1}^J \lambda(w_j,w_i) i(w_i,t|w_0,0).
 \label{eqn:cond_prob_flux}
 \end{equation}

Now we can define $\rho$ through an evolution law as follows
\begin{equation}
 \frac{d \rho(w_j,t|w_0,0)}{d t} =  \sum_{i=1}^J \lambda(w_j,w_i) i(w_i,t|w_0,0) - i(w_j,t|w_0,0) -\beta(w_j,t) \rho(w_j,t|w_0,0) .
\label{eqn:rho_evolution}
\end{equation}

Following Fedotov \cite{Fedotov2010} we can find an expression for the flux $i$ in terms of $\rho$ using Laplace transform methods on Eq. (\ref{eqn:rho_defn}) and Eq. (\ref{eqn:flux_def}) respectively. We first divide both equations by $e^{-\int\limits_0^{t} \beta(w_j,t'')dt''}$, this yields;
\begin{equation}
 \mathcal{L} \left\{ \rho(w_j,t|w_0,0) e^{\int\limits_0^t \beta(w_j,t'')dt''}\right\} = \mathcal{L} \left\{q(w_j,t|w_0,0) e^{\int\limits_0^t \beta(w_j,t')dt'}\right\} \mathcal{L} \{\phi(w_j,t) \}
\label{eqn:Laplace_rho}
\end{equation}
and
\begin{equation}
 \mathcal{L} \left\{ i(w_j,t|w_0,0) e^{\int\limits_0^t \beta(w_j,t'')dt''} \right\} = \mathcal{L} \left\{q(w_j,t|w_0,0) e^{\int\limits_0^t \beta(w_j,t')dt'} \right\} \mathcal{L} \{\psi(w_j,t) \}.
\label{eqn:Laplace_i}
\end{equation}
Rearranging Eqs. (\ref{eqn:Laplace_rho}) and (\ref{eqn:Laplace_i});
\begin{equation}
 \mathcal{L} \left\{i(w_j,t|w_0,0)e^{\int\limits_0^t \beta(w_j,t')dt'} \right\} = \mathcal{L} \left\{\rho(w_j,t|w_0,0) e^{\int\limits_0^t \beta(w_j,t')dt'} \right\} \frac{ \mathcal{L} \{\psi(w_j,t) \}}{ \mathcal{L} \{\phi(w_j,t) \}}.
\label{eqn:Laplace_i_rho}
\end{equation}
Inverting the Laplace transform, we get
\begin{equation}
 i(w_j,t|w_0,0) = \int\limits_0^t K(w_j,t-t') \rho(w_j,t'|w_0,0) e^{-\int\limits_{t'}^t \beta(w_j,t'') dt''} dt'
\label{eqn:inv_Laplace_i}
\end{equation}
where the memory kernel is defined by;
\begin{equation}
 K(w_j,t) = \mathcal{L}^{-1} \left\{ \frac{\mathcal{L}\{\psi(w_j,t) \}}{\mathcal{L} \{\phi(w_j,t) \}} \right\}.
\label{eqn:mem_ker}
\end{equation}
The master equation for the CTRW process on a network is found by the substituting Eq. (\ref{eqn:inv_Laplace_i}) into Eq. (\ref{eqn:rho_evolution})
\begin{eqnarray}
 \frac{d\rho(w_j,t|w_0,0)}{d t}  &=&  \sum_{i=1}^J \lambda(w_j,w_i) \int\limits_0^t K(w_i,t-t') \rho(w_i,t'|w_0,0) e^{-\int\limits_{t'}^t \beta(w_i,t'') dt''} dt'  \nonumber \\
 & &- \int\limits_0^t K(w_j,t-t') \rho(w_j,t'|w_0,0) e^{-\int\limits_{t'}^t \beta(w_j,t'') dt''}dt' \nonumber \\
 & & -\ \beta(w_j,t)\rho(w_j,t|w_0,t_0).
\label{eqn:rho_master_mid}
\end{eqnarray}
In the case of power law waiting time densities on a uniform grid network, this is similar to Eq. (29) in \cite{AYL2010}.

\subsection{Ensemble CTRW Birth-Death Master Equations}

To describe a birth-death process we also need to account for the creation of new particles, and hence need to consider ensembles of particles. We define $\eta(w_j,t)$ as the probability of a particle being created at vertex $w_j$ and at time $t$. Then we can define
\begin{equation}
u(w_j,t) = \sum\limits_{w_0\in W} \int\limits_0^t \rho(w_j,t|w_0,t_0) \eta(w_0,t_0) dt_0
\label{eqn:particle_dens}
\end{equation}
as the density of particles at vertex $w_i$  at time $t$.

Taking care to differentiate Eq. (\ref{eqn:particle_dens}) using Leibniz rule, we then substitute in Eq. (\ref{eqn:rho_master_mid}) and simplify to get
\begin{eqnarray}
\frac{d u(w_j,t)}{d t} & = & \sum\limits_{w_0\in W} \Biggl[ \eta(w_0,t) \rho(w_j,t|w_0,t) + \int\limits_0^t \eta(w_0,t_0) \frac{d \rho(w_j,t|w_0,t_0)}{d t}dt_0\Biggr]\\
		& =& \sum\limits_{w_0\in W} \int\limits_0^t \eta(w_0,t_0)  \Biggl[ \sum\limits_{i=1}^J \int\limits_0^t  K(w_i,t-t')\lambda(w_j,w_i) e^{-\int\limits_{t'}^t \beta(w_i,t'') dt''} \rho(w_i,t'|w_0,t_0)dt' \nonumber \\
		& & - \int\limits_0^t K(w_j,t-t') \rho(w_j,t'|w_0,0) e^{-\int\limits_{t'}^t \beta(w_j,t'') dt''}dt'  - \beta(w_j,t)\rho(w_j,t'|w_0,t_0)  \Biggr]dt_0 \nonumber \\
		& & +\ \eta(w_j,t).
\end{eqnarray}
As $\rho(w,t|w_0,t_0) =0 $ for all $t < t_0$ we can write;
\begin{equation}
		\begin{split} 
\frac{d u(w_j,t)}{d t} & = \int\limits_0^t \Biggl[ \sum\limits_{i=1}^J K(w_i,t-t') \lambda(w_j,w_i) e^{-\int\limits_{t'}^t \beta(w_i,t'')dt''}  \sum\limits_{w_0\in W} \int_0^{t'} \rho(w_i,t'|w_0,t_0)\eta(w_0,t_0)dt_0  \\
		& \quad - e^{-\int\limits_{t'}^t \beta(w_j,t'')dt''} K(w_j,t-t') \sum\limits_{w_0 \in W} \int\limits_0^{t'} \rho(w_j,t'|w_0,t_0)\eta(w_0,t_0)dt_0 \Biggr] dt'  \\
		& \quad - \beta(w_j,t) \sum\limits_{w_0 \in W} \int\limits_0^{t} \rho(w_j,t|w_0,t_0) \eta(w_0,t_0) dt_0 + \eta(w_j,t).  \end{split}
\label{eqn:u_gen_master}
\end{equation}
Finally we arrive at the master equation for a CTRW with reactions on networks
\begin{eqnarray}
 \frac{d u(w_j,t)}{d t} & =& \int\limits_0^t  \Biggl[ \sum\limits_{i=1}^J K(w_i,t-t') \lambda(w_j,w_i) e^{-\int\limits_{t'}^t \beta(w_i,t'')dt''} u(w_i,t')  \nonumber  \\ 
		& &  - K(w_j,t-t') e^{-\int\limits_{t'}^t \beta(w_j,t'')dt''} u(w_j,t') \Biggr] dt' - \beta(w_j,t)u(w_j,t) +  \eta(w_j,t)
\label{eqn:master_CTRW}
\end{eqnarray}
which can be written in the general form of a reaction-diffusion equation
\begin{equation}
\frac{d u(w_j,t)}{d t}  = L[u(w_j,t)] + f(u(w_j,t))
\label{eqn:gen_master_CTRW_0}
\end{equation}
where 
 \begin{eqnarray}
  L [u(w_j,t)] & =& \int\limits_0^t  \Biggl[ \sum\limits_{i=1}^J K(w_i,t-t') \lambda(w_j,w_i) e^{-\int\limits_{t'}^t \beta(w_i,t'')dt''} u(w_i,t')  \nonumber \\ 
		& & \qquad - \ K(w_j,t-t') e^{-\int\limits_{t'}^t \beta(w_j,t'')dt''} u(w_j,t') \Biggr] dt'
\label{eqn:gen_master_CTRW}
\end{eqnarray}
is the CTRW network Laplacian and 
\begin{equation}
f(u(w_j,t)) = -\beta(w_j,t)u(w_j,t) +  \eta(w_j,t).
\end{equation}
models the reaction kinetics on a vertex expressed in terms of birth and death processes. Previously it has been noted that the reaction kinetics are incorporated into the transport operator for a CTRW process \cite{HLW2006}. It can be seen from Eq. (\ref{eqn:gen_master_CTRW}), that here only the death processes are incorporated in the transport operator. 

\subsection{CTRW Laplacian with Exponential Waiting Times}
\label{sec:laplace}

We now apply Eq. (\ref{eqn:master_CTRW}) to the case of exponential waiting times,
\begin{equation}
 \psi(w_j,t) = \alpha(w_j) e^{-\alpha(w_j) t}.
\label{eqn:exp_wait_time}
\end{equation}
This greatly simplifies the master equation. In general, the Laplace transform of the waiting time density, $\bar{\psi} (w_j, s) = \mathcal{L} \{ \psi(w_j, t)\}$ and the Laplace transform of the survival probability $\bar{\phi} (w_j, s) = \mathcal{L} \{ \phi(w_j, t)\}$ are related by $\bar{\phi}(w_j, s) = \frac{1-\bar{\psi}(w_j,s)}{s}$.

As $\psi(w_j,s)$ is exponential, then $\bar{\psi}(w_j, s) = \frac{s}{s+\alpha(w_j)}$ so $\frac{\bar{\psi}(w_j,s)}{\bar{\phi}(w_j,s)} = \alpha(w_j)$ with the inverse Laplace transform $\alpha(w_j) \delta(t)$ and thus 
\begin{equation}
 K(w_j,t) = \alpha(w_j)\delta(t).
\label{eqn:k_explic}
\end{equation}
By substitution, we can rewrite the master equation, Eq. (\ref{eqn:master_CTRW}), as 
\begin{align}
\frac{d u(w_j,t)}{d t} & = \sum_{i=1}^J \alpha(w_i) \lambda(w_j,w_i)  u(w_i,t) 
   - \alpha(w_j)  u(w_j,t) + f(u(w_j,t)) \\
  & = \sum\limits_{i=1}^J L_{i,j} u(w_i,t) +f(u(w_j,t))
\label{eqn:rho_master_fin}
\end{align}
where $L$ is a member of the family of general CTRW network Laplacians defined as
\begin{equation}
L_{i,j}= \left\{ \alpha(w_i)\lambda(w_j,w_i) - \alpha(w_j) \delta_{i,j} \right\}.
\label{eqn:CTRW_lap}
\end{equation}

In the following we consider two special cases of Eq. (\ref{eqn:CTRW_lap}). For both cases, we assume that the jump probability and the waiting time probability are only functions of vertex degree and time respectively. Both of these cases have been previously considered without reactions \cite{SDM2008}. To describe the network we used the adjacency matrix \begin{align}
A_{i,j}= \left\{ \begin{array}{c c}1 & \mathrm{if \ vertices\ } w_i \mathrm{ \ and \ } w_j \mathrm{\ are \ connected}\\
0 & \mathrm{otherwise.}
\end{array}\right.
\end{align}
Each vertex, $w_j$, has a degree, $k_{j}$, that is defined as the total number of edges that link to the vertex. This is also the sum of each row of the adjacency matrix.

\subsubsection{Case A}
We assume that the waiting time on each vertex is identical and that the probability of jumping across an edge is equal for a given vertex. Formally, we let $\alpha(w_j) =  \alpha_A$ for all $w_j \in W$ and let $\lambda(w_j,w_i) = \frac{1}{k_i}$. Thus our Laplacian becomes 
\begin{equation}
L_{i,j} = \alpha_{A} \left( \frac{1}{k_i} A_{i,j} - \delta_{i,j} \right).
\label{eqn:Case_a_lap}
\end{equation}
Note that in this case the Laplacian is not symmetric. In general, the steady state for a CTRW with no reactions using this Laplacian will not be uniform on the vertex set.

\subsubsection{Case B}
An alternative to Case A, is to let the waiting time on each vertex change proportionally to the vertex degree. This allows the rate of particles jumping along each edge to be constant. Formally, we let $\alpha(w_j) =  \alpha_B k_j$ for all $w_j \in W$ and let $\lambda(w_j,w_i) = \frac{1}{k_i}$. Thus our Laplacian can be described as 
\begin{equation}
L_{i,j} = \alpha_B \left( A_{i,j} - k_j \delta_{i,j} \right).
\label{eqn:Case_b_lap}
\end{equation}
This is the well studied graph Laplacian \cite{Nakao2010,Hata2012}. It is important to note that if a connected network is regular, so that $k_j$ is the same for all $w_j$, then the cases are equivalent up to a scale factor.

\section{Pattern Formation}
\label{sec:patterns}

In continuum reaction-diffusion partial differential equations the concentration may vary on the spatial domain. This patterning also holds on discretization of the spatial manifold where the Laplacian operator is replaced with a discrete Laplace Beltrami operator \cite{Xu2004} \cite{PSPBM2004}. This can be considered as a special case of pattern formation on a discrete network. In general variation in concentrations on each vertex may occur on any network. This variation may arise in a number of different ways. First, if the Laplacian matrix is not symmetric in a system without reactions, there will be a buildup of concentrations on vertices according to their degree. This steady state pattern will be a multiple of the eigenvector of the Laplacian with zero eigenvalue. Second, with a nonsymmetric Laplacian matrix in systems where the reactions have a finite nonzero steady state solution, the interplay between the diffusion and the reactions will cause a different pattern across the network.  In this pattern, unlike the no reaction case, vertices with the same degree may have different concentrations. We refer to these two mechanisms as Laplacian pattern formation, as the patterning is driven by the nonsymmetric Laplacian matrix. Last, if the reaction terms permit a Turing instability, whereby the spatially homogeneous steady state becomes unstable in the presence of diffusion, then a Turing pattern may form \cite{Turing1952}. This instability is permitted for both a symmetric and a nonsymmetric Laplacian matrix. 

\subsection{Laplacian Pattern Formation}

To completely eliminate any interplay with Turing patterns, we consider a single species model that cannot permit a Turing instability \cite{Murray2002}.
\begin{equation}
\frac{d \mathbf{u}}{d t}  = \mathbf{L}. \mathbf{u} + \mathbf{f}(\mathbf{u})\\
\label{eqn:single_species}
\end{equation}
where $\mathbf{u} = \left( \begin{array}{ccc}
u(w_1,t), \ \hdots,\ u(w_J,t) \end{array}\right)^T$ and $\mathbf{f}(\mathbf{u}) = \left( \begin{array}{ccc}
f(u(w_1,t)),\ \hdots ,\ f(u(w_J,t)) \end{array}\right)^T$.

If we take the trivial reaction term, $\mathbf{f}(\mathbf{u})=\mathbf{0}$, the only possible form of spatial pattern formation comes from a nonsymmetric Laplacian. The master equation is then simply
\begin{equation}
\frac{d \mathbf{u}}{d t}  = \mathbf{L}. \mathbf{u} .
\label{eqn:single_species_no_f}
\end{equation}
It is clear that if $L$ is nonsymmetric, there will be a nonuniform rate of particle transport along the network that produces a pattern of different concentrations of particles in the long term. This concentration vector must be a multiple of the eigenvector of the Laplacian matrix with eigenvalue zero as it is a solution of the the vector $\mathbf{u}$ that satisfies the steady state of Eq. (\ref{eqn:single_species_no_f}), i.e.,
\begin{equation}
\mathbf{L} .\mathbf{u} = \mathbf{0}.
\label{eqn:great_evec}
\end{equation}
If the reaction term $\mathbf{f}(\mathbf{u}) \neq \mathbf{0}$ and has a finite non-zero equilibrium solution then the pattern will no longer correspond to an eigenvector of the Laplacian matrix. 

Define the reaction steady state $\mathbf{u}^*$ such that the reaction term $\mathbf{f}(\mathbf{u}^*)=\mathbf{0}$. By considering the behaviour of $\mathbf{u} = \mathbf{u}^* + \Delta \mathbf{u}$ where $\Delta \mathbf{u}$ is some perturbation, Eq. (\ref{eqn:single_species}) may be linearized to become 
\begin{equation}
\frac{d \Delta\mathbf{u}}{d t}  =\mathbf{D}\mathbf{f}(\mathbf{u}^*). \Delta\mathbf{u}+ \mathbf{L}.\mathbf{u}^* + \mathbf{L}.\Delta\mathbf{u}.
\label{eqn:single_species_per}
\end{equation}
The steady state solutions are given by 
\begin{equation}
\label{eqn:lp}
 \Delta\mathbf{u} =-\left( \mathbf{L}+\mathbf{D}\mathbf{f}(\mathbf{u}^*)\right)^{-1}.\mathbf{L}. \mathbf{u}^{*},
 \end{equation}
 provided the inverse exists, where $\mathbf{D}\mathbf{f}(\mathbf{u}^*)$ is the derivative operator of the reaction vector field $\mathbf{f}$ at $\mathbf{u}^*$. The perturbed state $\mathbf{u}^* + \Delta \mathbf{u}$ is the linear prediction of the pattern. In the case where the reaction term is linear in $\mathbf{u}$ this is the exact solution for the pattern. It should also be noted that in the case of a symmetric Laplacian, $\mathbf{L} .\mathbf{u}^{*}=\mathbf{0}$, as the constant vector is always a null vector of a symmetric Laplacian. This equation does not hold without reactions as a nonsymmetric Laplacian is in general non-invertible.

This type of pattern formation is very different from a Turing pattern formation. A Turing pattern is formed by an instability in the dynamics whereby an otherwise stable solution is made unstable by the presence of diffusion. It is a bifurcation phenomena as the diffusion must reach a critical value before the solution becomes unstable. In this case however we have a pattern that will form for any amount of diffusion, the stability of the solution does not change but rather the solution itself is a function of the diffusion. To examine this we introduce a scale parameter $s$ for our Laplacian whereby the linear predictor for the pattern becomes
\begin{equation}
\Delta\mathbf{u} =-\left( s \mathbf{L}+\mathbf{D}\mathbf{f}(\mathbf{u}^*)\right)^{-1}s \mathbf{L} .\mathbf{u}^{*}
\label{eqn:lin_Du}
\end{equation}
The new scalar parameter governs the speed of diffusion in the system. 
In the limit $s\rightarrow0$ we have;
\begin{equation}
\Delta\mathbf{u} =\mathbf{0}
\end{equation}
which corresponds to no pattern, and each vertex taking the equilibrium values of the reactions.
In the limit $s\rightarrow \infty$ care has to be taken with the existence of the inverse. Taking Eq. (\ref{eqn:lin_Du}), and multiplying through by $\left( s \mathbf{L}+\mathbf{D}\mathbf{f}(\mathbf{u}^*)\right)$, dividing by $s$, and taking the limit of $s \rightarrow \infty$, we have
\begin{align}
\lim_{s \rightarrow \infty}\left( \mathbf{L}+\frac{\mathbf{D}\mathbf{f}(\mathbf{u}^*)}{s}\right)\Delta\mathbf{u} &=- \mathbf{L} .\mathbf{u}^{*},\\
\mathbf{L}.\Delta\mathbf{u}+\mathbf{L} .\mathbf{u}^{*}&=\mathbf{0},\\
\mathbf{L}.\mathbf{u}&=\mathbf{0}
\end{align}
and so the pattern will be equivalent to the case with no reactions. In this way it can be seen that this type of pattern formation is the mixing of the pattern formation from the nonsymmetric Laplacian and the reaction equilibrium.

\subsection{Turing Pattern Formation}
We consider a general two species reaction-diffusion system with particle concentrations $u$ and $v$. The master equations, [see Eq. (\ref{eqn:rho_master_fin})], can be written as
\begin{equation}
\frac{d  u(w_j,t)}{d t} = \sum\limits_{i=1}^J \alpha_u L_{i,j} u(w_i,t) + f(u(w_j,t),v(w_j,t))
\label{eqn:master_u}
\end{equation}
and
\begin{equation}
\frac{d v(w_j,t)}{d t} = \sum\limits_{i=1}^J \alpha_v L_{i,j} v(w_i,t) + g(u(w_j,t),v(w_j,t))
\label{eqn:master_v}
\end{equation}
where the functions $f(u(w_j,t),v(w_j,t))$ and $g(u(w_j,t),v(w_j,t))$ incorporate the creation and destruction probabilities, i.e.,
\begin{align}
f(u(w_j,t),v(w_j,t)) & = \eta_u(w_j,t) - \beta_u(w_j,t)u(w_j,t) \\
			    & = \eta_u(u(w_j,t),v(w_j,t)) - \beta_u(u(w_j,t),v(w_j,t))u(w_j,t)
\label{eqn:u_expansion}
\end{align}
and similarly for $g(u(w_j,t),v(w_j,t))$. Here note that the $\alpha_u$ and $\alpha_v$ are factored out of the Laplacian, $L_{i,j}$, so that the operator is the same in both equations. 

\subsubsection{Linear Stability Analysis}

To consider Turing instabilities, we first rewrite Eqs. (\ref{eqn:master_u}) and (\ref{eqn:master_v}) in vector form 
\begin{equation}
\frac{d \mathbf{X}}{d t} = \mathbf{\Lambda }\mathbf{X} + \mathbf{F} (\mathbf{X})
\label{eqn:vec_gen}
\end{equation}
where 
\begin{align*}
\mathbf{X} & = \left( \begin{array}{cccccc}
u(w_1,t) ,\hdots, u(w_J,t), \ v(w_1,t), \hdots, v(w_J,t) \end{array}\right)^T \\
& = \left( \begin{array}{cccccc} X_1, \hdots, X_J, \ X_{J+1}, \hdots, X_{2J} \end{array}\right)^T,\\
\mathbf{F}(\mathbf{X}) &= \left( \begin{array}{cccccc} f(X_1,X_{J+1}),\hdots ,f(X_J,X_{2J} ),\ 		g(X_1,X_{J+1}), \hdots ,g(X_J,X_{2J}) \end{array}\right)^T \\
	& = \left( \begin{array}{cccccc} F_1, \hdots ,F_J, \ F_{J+1}, \hdots, F_{2J} \end{array}\right)^T ,\\
\textrm{and \ } \mathbf{\Lambda} &= \left( \begin{array}{cc} \alpha_u \mathbf{L}\ \  \mathbf{0} \ \\ \mathbf{0} \ \ \ \alpha_v \mathbf{L} \end{array} \right).
\end{align*}

Linearising about the steady state $\mathbf{X}^*$ with $\mathbf{F}(\mathbf{X}^*)=\mathbf{0}$ and $\mathbf{X} = \mathbf{X}^* + \Delta \mathbf{X}$. We then substitute this into Eq. (\ref{eqn:vec_gen}) to get 
\begin{equation}
\frac{d \Delta \mathbf{X}}{d t}  = \left( \mathbf{\Lambda} + \mathbf{DF}(\mathbf{X}^*)\right). \mathbf{\Delta}\mathbf{X} + \mathbf{\Lambda}. \mathbf{X}
\label{eqn:lin_gen}
\end{equation}
where 
\begin{equation} 
\mathbf{D}\mathbf{F}(\mathbf{X}^*)_{i,j} = \left. \frac{d F_i}{d X_j}\right| _{\mathbf{X}^*}.
\end{equation}
We now apply the affine transform $\Delta \mathbf{Y} = \Delta \mathbf{X} + \left( \mathbf{\Lambda} + \mathbf{DF}(\mathbf{X}^*)\right)^{-1} .\mathbf{\Lambda}.\mathbf{X}^*$ to Eq. (\ref{eqn:lin_gen}) yielding
\begin{equation}
\frac{d \Delta \mathbf{Y}}{d t} = \left( \mathbf{\Lambda} +\mathbf{DF}(\mathbf{X}^*) \right). \Delta \mathbf{Y} 
\end{equation}
which gives solutions of the form 
\begin{equation}
\Delta \mathbf{Y} = \sum\limits_{j=1}^J e^{\mu_j t} P_j(t) \boldsymbol{\nu}_j ,
\label{eqn:delta_y_soln}
\end{equation}
where $\mu_j$ is the $j^{th}$ eigenvalue of $\left( \mathbf{\Lambda}+\mathbf{DF}(\mathbf{X}^*)\right)$ with corresponding eigenvector $\boldsymbol{\nu}_j$ and $P_j$ is a polynomial in $t$ for repeated eigenvalues. The long time behaviour of $\Delta \mathbf{Y}$ is then approximated by 
\begin{equation}
\Delta \mathbf{Y}  \sim e^{\mu^* t} P^*(t) \boldsymbol{\nu}^* 
\label{eqn:delta_y_app}
\end{equation}
where $\mu^*$ corresponds to the eigenvalue with the largest real component. 

In the linear stability analysis the concentrations of the two species evolve as 
\begin{equation}
\mathbf{X} \sim e^{\mu^* t} P(t) \boldsymbol{\nu}^* + \mathbf{X}^* -  \left( \mathbf{\Lambda} + \mathbf{DF} (\mathbf{X}^*)\right)^{-1}.\mathbf{\Lambda}.\mathbf{X}^*
\label{eqn:lin_stab_gen}
\end{equation}
where $ \left(\mathbf{\Lambda} + \mathbf{DF} (\mathbf{X}^*)\right)^{-1}.\mathbf{\Lambda}.\mathbf{X}^*=\mathbf{0}$ if the Laplacian is symmetric.

In the continuum case the real components of the eigenvalues (i.e., stability) of the homogeneous steady state can be plotted against spatial frequency to obtain a dispersion relation showing the range of spatial frequencies that will grow with time. In an analogous manner for the network case, we plot the stability of the reaction-diffusion system as a function of a scale parameter $s$ equivalent to scaling the waiting time for both species on the network:
\begin{equation}
\frac{d \mathbf{X}}{d t} =s\mathbf{ \Lambda}. \mathbf{X} + \mathbf{F} (\mathbf{X})
\label{eqn:vec_gen_s}
\end{equation}
When $s$ is zero there is no coupling between vertices and each vertex will be in equilibrium according to the reaction equations. We identify two critical values of $s$ greater than zero from our linear stability analysis. First, the Laplacian type pattern arising from the coupling of reactions to diffusion may change to a Turing pattern at a critical value of $s$.  As $s$ is further increased the Turing pattern, if it occurs, will persist untill the second critical value when the pattern reverts back to a Laplacian type pattern.

\section{Examples of Pattern Formation}
\label{sec:numerical}

To illustrate the properties of both the Laplacian and the Turing patterns on networks the master equations with exponential waiting times were solved numerically. For the Laplacian patterns, Eq. (\ref{eqn:rho_master_fin}) was solved with logistic reaction kinetics using the Case A Laplacian operator. The Turing patterns were examined by solving Eqs. (\ref{eqn:master_u}) and (\ref{eqn:master_v}) with Gierer-Meinhardt reaction kinetics and both the Case A and the Case B Laplacian operators. In both cases the equations were solved on random networks generated by the Barab\'asi-Albert algorithm \cite{Barabasi1999}, and the Watts-Strogatz algorithm \cite{Watts1998}.

The first reaction kinetics that we consider is the logistic equation \cite{Verhulst1845}, 
\begin{equation}
f\left(u(x,t)\right)= r u(x,t) (1-u(x,t))
\label{eqn:logistic}
\end{equation}
where $r$ is a constant. This governs the growth rate of a single species. In the following examples we take $r=1$. When applied to a network, this simple example could be considered a model for animal populations in a set of connected habitats, where the population is constrained by natural limits. 
The reaction-diffusion master equation in this case is found by substituting Eq. (\ref{eqn:logistic}) into Eq. (\ref{eqn:single_species}),
\begin{equation}
\frac{d u(w_j,t)}{d t} = \sum_{i=1}^J L_{i,j} u(w_i,t) + u(w_i,t) (1-u(w_i,t))
\end{equation}
Here $L$ is the Case A Laplacian with $\alpha=1$, i.e., $L_{i,j}=\frac{A_{i,j}}{k_{i}}-\delta_{i,j}$. 

The second model we consider has Gierer-Meinhardt reaction kinetics \cite{GM1972}. This is a two species model that permits Turing instabilities. The reaction terms in the model are
\begin{equation}
f(u(x,t),v(x,t))=c\rho \ \frac{u(x,t)^{2}}{v(x,t)} - \mu \ u(x,t) + \rho_{0} \rho 
\label{eqn:gm_f}
\end{equation}
and
\begin{equation}
g(u(x,t),v(x,t))=c_d \rho \ u(x,t)^{2} - \nu\  v(x,t).
\label{eqn:gm_g}
\end{equation}
The Gierer-Meinhardt reaction  kinetics are only valid for $u(x,t)$ and $v(x,t)>0$ for all $x$ and $t$. Furthermore we assume that all vertices have equal volume so that the number of particles and the concentration are interchangeable.

We can apply this model to a network setting by substituting Eqs. (\ref{eqn:gm_f}) and (\ref{eqn:gm_g}) in Eqs. (\ref{eqn:master_u}) and (\ref{eqn:master_v}):
\begin{align}
\frac{d u(w_j,t)}{d t}&=c \rho \ \frac{u(w_j,t)^{2}}{v(w_j,t)} - \mu \ u(w_j,t) + \rho_{0} \rho +\alpha_u\sum_{i=1}^J L_{i,j} u(w_i,t)  \label{eqn:u_gm},\\
\frac{d v(w_j,t)}{d t}&=c_d \rho\ u(w_j,t)^{2} - \nu \ v(w_j,t) +\alpha_v \sum_{i=1}^J L_{i,j} v(w_i,t) . 
\label{eqn:v_gm}
\end{align}
In the following we use
\begin{equation}
\rho_0 =1, \ \rho = 1, \ \nu=\frac{7}{32},\ \mu=\frac{5}{256},\ c=1, \textrm{and} \ c_d=\frac{5}{128}.
\end{equation}
For the Case A Laplacian we use
\begin{equation}
\alpha_u=1 \ \mathrm{and} \ \alpha_v=\frac{1}{256}.
\end{equation}
To place the Case B Laplacian, i.e., $L_{i,j}=A_{i,j}-k_{i}\delta_{i,j}$, on a comparable footing to Case A we rescale the parameters $\alpha_u$ and $\alpha_v$ to ensure that the mean waiting time across the network is comparable in both cases. Explicitly 
\begin{equation}
\alpha_u = \frac{J}{\sum_{j=1}^{J}k_{j} }\ \mathrm{and} \ \alpha_v =\left( \frac{1}{256}\right)\frac{J}{\sum_{j=1}^{J}k_{j} }.
\end{equation}

\subsection{Barab\'asi-Albert Network}

The Barab\'asi-Albert (BA) network is a random network with a power law distribution of vertex degrees \cite{Barabasi1999}. The network is iteratively generated by adding a vertex at each step that connects to $k$ existing vertices where the probability of attachment is proportional to the degree of the existing vertices. 

\begin{figure}[h]
\begin{center}
\includegraphics[width=1\linewidth]{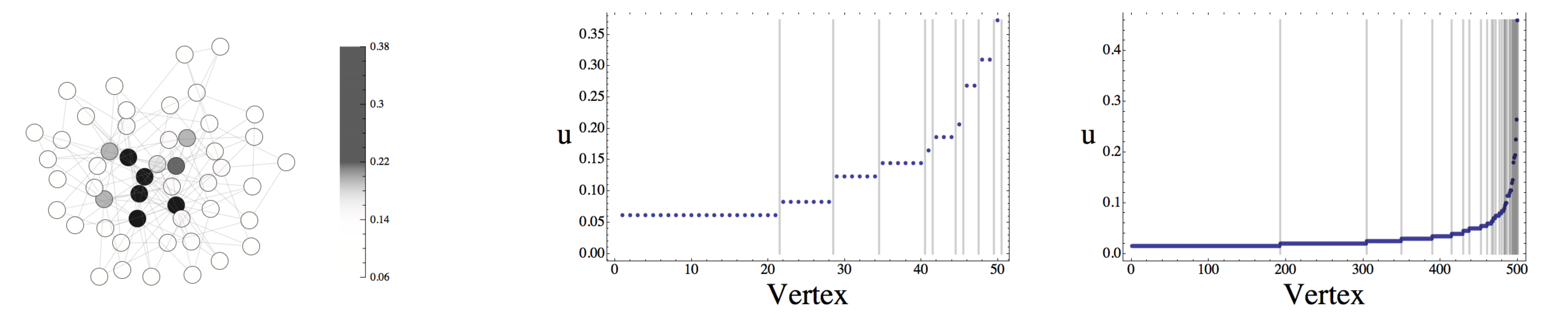}
\caption{Case A Laplacian diffusion on a BA network with $k=3$. A characteristic pattern of a network with 50 vertices (\emph{left}) and the distribution of concentrations on a 50 (\emph{centre}) and a 500 (\emph{right}) vertex network. The vertical lines demark regions of different vertex degree}
\label{fig:BA_500_3_Evec}
\end{center}
\end{figure}

For completeness, we first consider a purely diffusive process governed by the Case A Laplacian on BA networks. In this case,  the concentration on each vertex is proportional to its degree, resulting in the shape of the concentrations to be similar to a power law shape as shown in Fig. \ref{fig:BA_500_3_Evec}. This distribution, which corresponds to the eigenvector of the Laplacian with a zero eigenvalue, is monotonic with increasing vertex degree.

With the addition of the logistic reaction term, concentrations across the network change.
\begin{figure}[h]
\begin{center}
\includegraphics[width=1\linewidth]{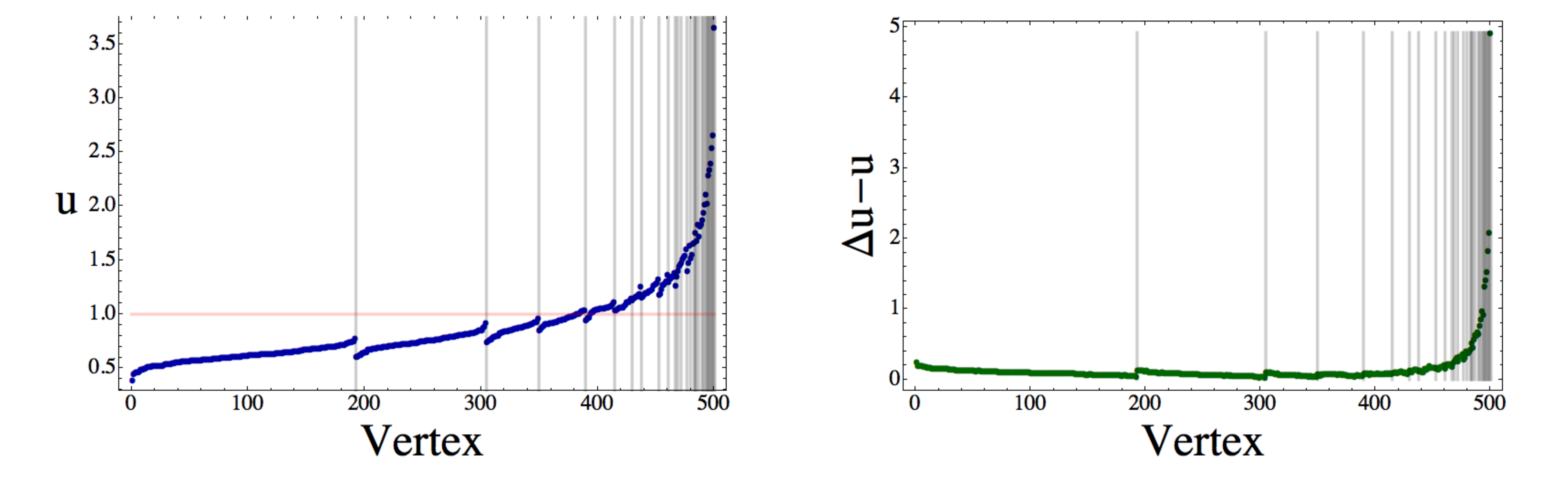}
\caption{ Case A Laplacian with logistic reaction kinetics on a BA network with $k=3$. Concentration (\emph{left}) and difference of linear predictor and concentration (\emph{right}). The pale horizontal line is the reaction steady state value and the data is ordered according to concentration within each segment of constant vertex degree. The vertical lines demark regions of different vertex degree}
\label{fig:BA_500_3_Log}
\end{center}
\end{figure}
The overall power law shape is still present, however the pattern is no longer monotonic with respect to vertex degree. The linear predictor of the pattern across the network, Eq. (\ref{eqn:lp}), is a good approximation near the reaction steady state, see Fig. \ref{fig:BA_500_3_Log}. The mean of the concentration across the network is reduced from the reaction equilibrium value.

\begin{figure}[h]
\begin{center}
\includegraphics[width=1\linewidth]{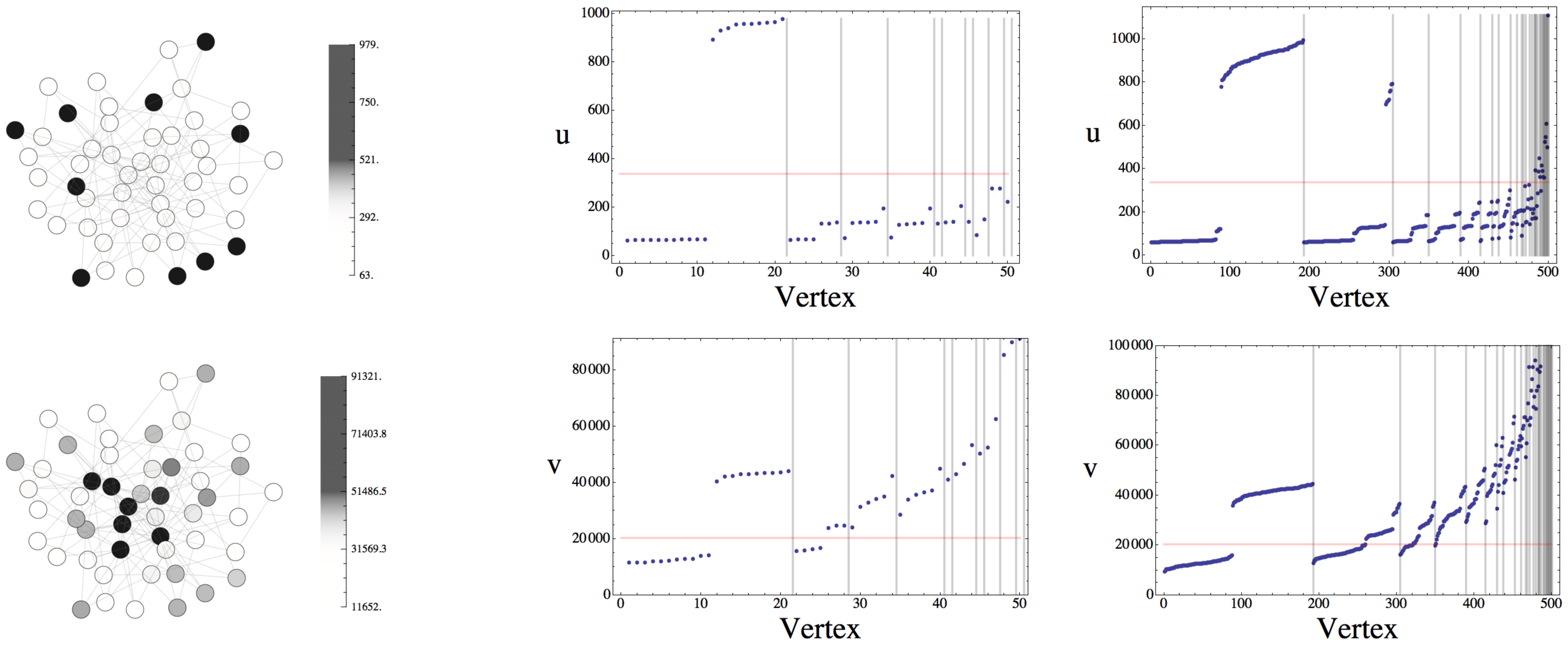}
\caption{ Case A Laplacian Gierer-Meinhardt reaction diffusion on a BA network with $k=3$. The characteristic pattern of $u$ (\emph{top}) and $v$ (\emph{bottom}) on 50 (\emph{left}) vertex network and the corresponding distribution of concentrations on a network with 50 vertices (\emph{centre}) and 500 (\emph{right}) vertices. The horizontal line gives the reaction steady state.}
\label{fig:BA_500_3_GM_CA}
\end{center}
\end{figure}

When the reactions are governed by the Gierer-Meinhardt model, Turing instabilities can arise producing patterns different to the Laplacian patterns. We first consider the patterns when the Case A Laplacian is used; a representative example is shown in Fig. \ref{fig:BA_500_3_GM_CA}.
The exact pattern is determined by the initial conditions of the system. This multi-stability is in contrast to the Laplacian patterns that have no initial condition dependence. In all observed Turing patterns the concentrations for both types of particles are split on vertices with low degree. There is also an increasing concentration as a function of vertex degree as seen in the previous BA patterns, especially for the concentration of $u$ particles. 

\begin{figure}[h]
\begin{center}
\includegraphics[width=1\linewidth]{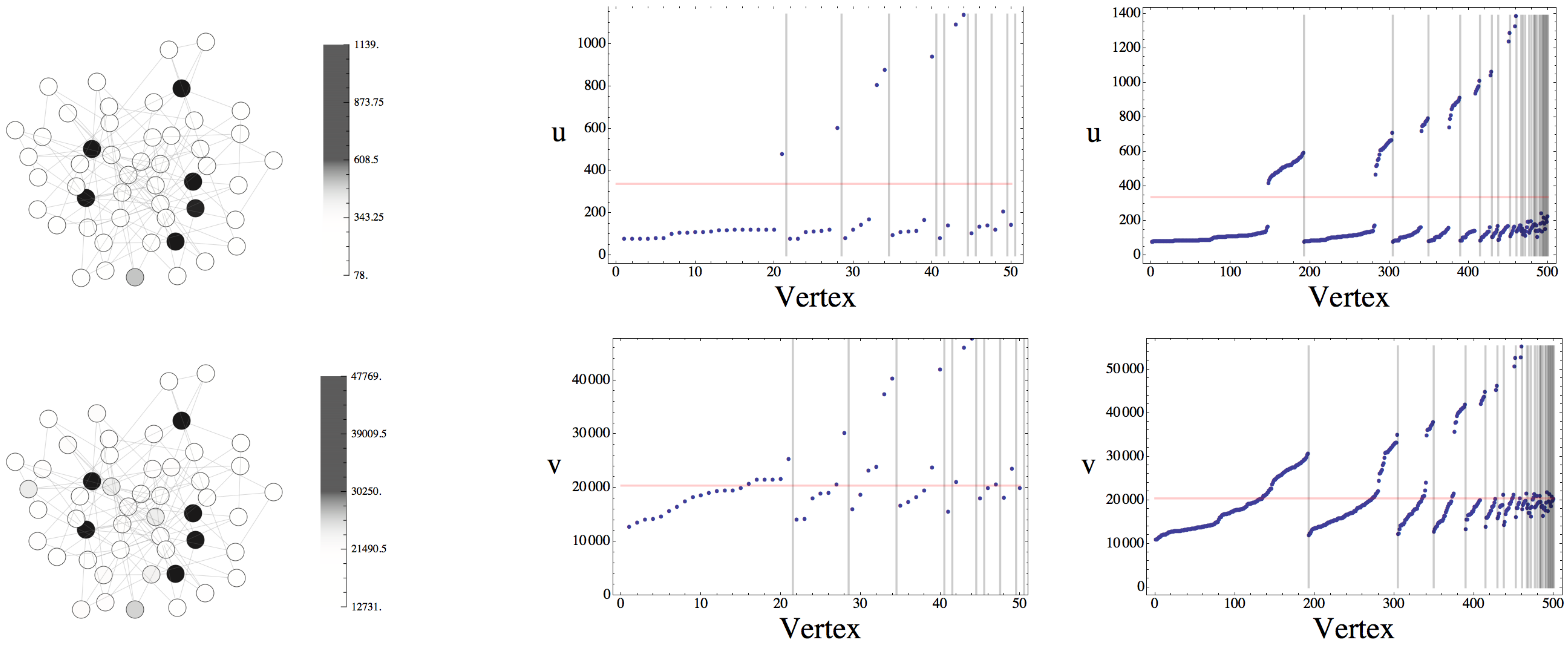}
\caption{ Case B Laplacian Gierer-Meinhardt reaction diffusion on a BA network with $k=3$. The characteristic pattern of $u$ (\emph{top}) and $v$ (\emph{bottom}) on 50 (\emph{left}) vertex network and the corresponding distribution of concentrations on a network with 50 vertices (\emph{centre}) and 500 (\emph{right}) vertices. The horizontal line gives the reaction steady state}
\label{fig:BA_500_3_CB}
\end{center}
\end{figure}

When considering the Case B Laplacian, as shown in Fig. \ref{fig:BA_500_3_CB}, the patterns shown are clearly very different to those in Fig. \ref{fig:BA_500_3_GM_CA}. There is a splitting of concentration across all vertex degree segments and the lower level shows practically no increase in concentration as a function of vertex degree. Moreover, the concentrations achieved are higher and lower for $u$ and $v$, respectively, when compared to Case A. Unlike for the Case A model, the patterns exhibited by $u$ and $v$ have similar profiles. 

It is interesting to note the change in mean concentration for the two Turing patterns. For the Case A pattern we have both the the mean of $u$ and $v$ greater then the reaction steady state value with no diffusion. For Case B, the opposite is true and the mean values are less then the reaction steady state values. 

\begin{figure}[h]
\begin{center}
\includegraphics[width=1\linewidth]{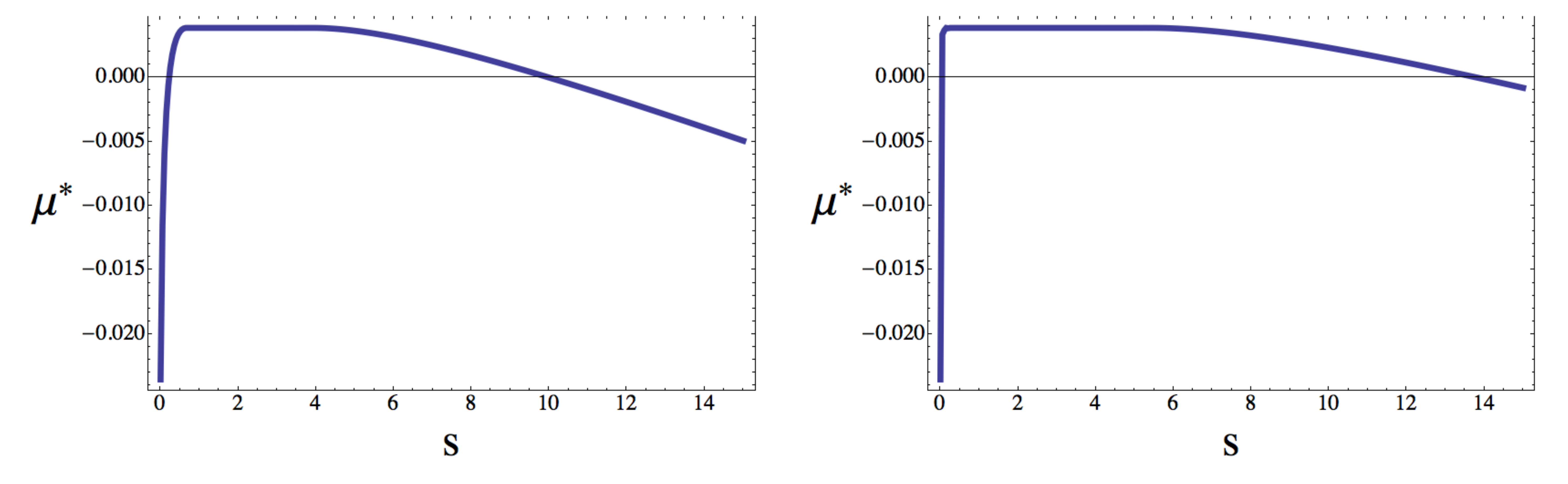}
\caption{ Dispersion relation for Case A (\emph{left}) and Case B (\emph{right}) Laplacian Gierer-Meinhardt reaction diffusion on a 500 vertex BA network with $k=3$.}
\label{fig:BA_500_3_GM_D}
\end{center}
\end{figure}

The dispersion relations for Cases A and B are plotted in Fig. \ref{fig:BA_500_3_GM_D}.
This figure shows that the Turing instability occurs over a larger range of the scale parameter $s$ in Case B than in Case A. 

\subsection{Watts-Strogatz Network}

The Watts-Strogatz (WS) network is a random network characterized by the small world property of having a low graph diameter; the length of the longest path between any two vertices that does not involve loops or back tracking \cite{Watts1998}. The network is generated by creating a ring lattice where each vertex is connected to $k$ adjacent vertices. Then each subsequent edge may be reconnected with some probability $p$ to some other vertex chosen uniformly from all vertices. 

\begin{figure}[h]
\begin{center}
\includegraphics[width=1\linewidth]{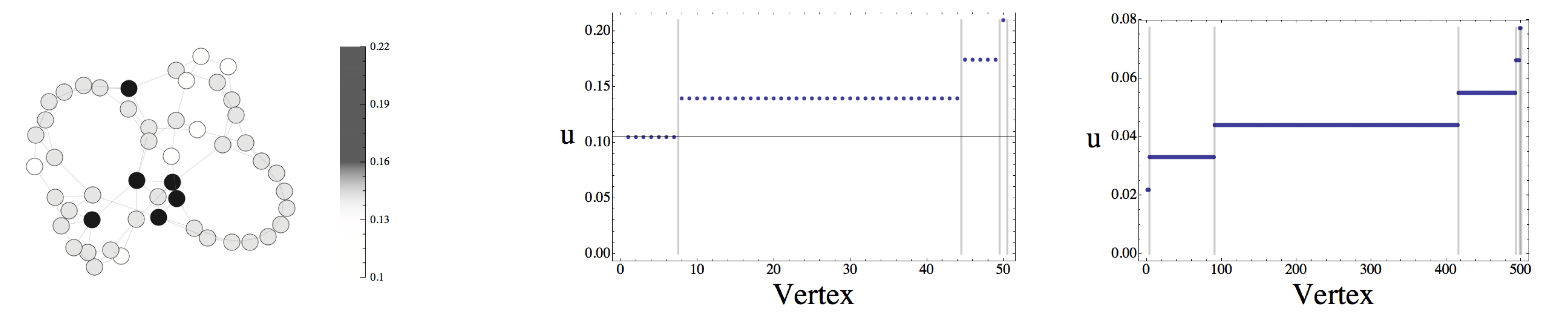}
\caption{Case A Laplacian diffusion on a WS network with $k=3$ and $p=0.1$. A characteristic pattern of a network with 50 vertices (\emph{left}) and the distribution of concentrations on a 50 (\emph{centre}) and a 500 (\emph{right}) vertex network. The vertical lines demark regions of different vertex degree.}
\label{fig:WS_500_Evec}
\end{center}
\end{figure}

Once again we first consider a purely diffusive process governed by the Case A Laplacian. The concentration on each vertex is proportional to the degree of the vertex as shown in Fig. \ref{fig:WS_500_Evec}.

\begin{figure}[h]
\begin{center}
\includegraphics[width=1\linewidth]{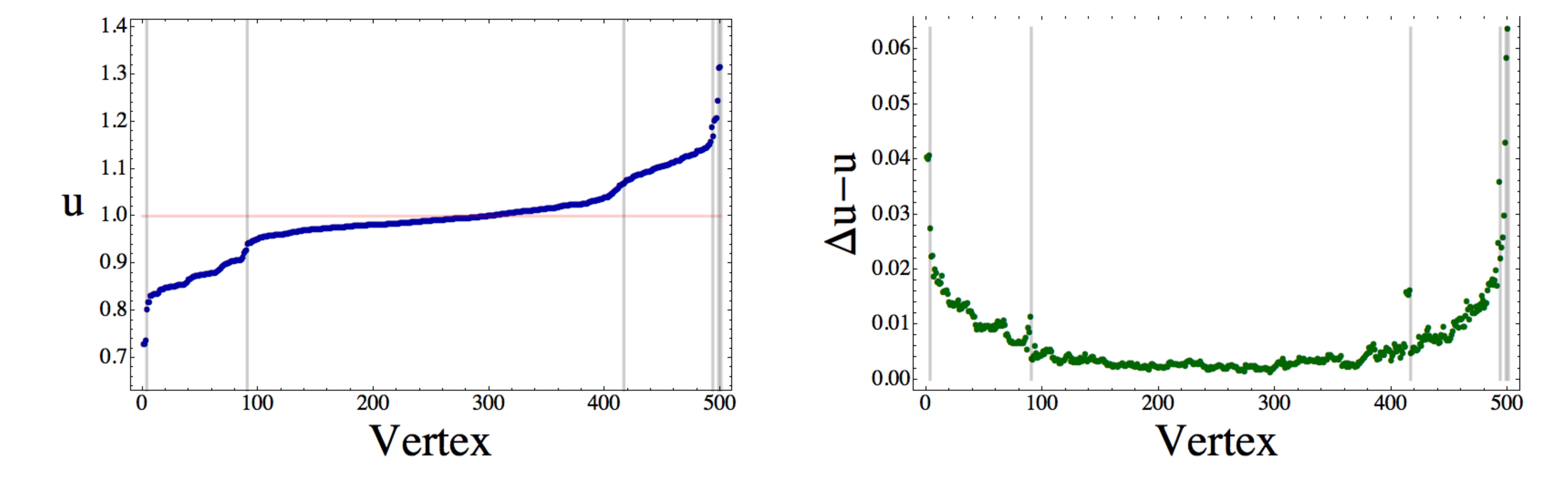}
\caption{ Case A Laplacian with logistic reaction kinetics on a WS network with $k=3$ and $p=0.1$. Concentration (\emph{left}) and difference of linear predictor and concentration (\emph{right}). The pale horizontal line is the reaction steady state value and the data is ordered according to concentration within each segment of constant vertex degree. The vertical lines demark regions of different vertex degree.}
\label{fig:WS_500_Log}
\end{center}
\end{figure}

The same network with logistic reaction kinetics is shown in Fig. \ref{fig:WS_500_Log}.
Similarly to the BA network, each segment of vertices with identical degree shows a gentle slope in concentration with tapered boundaries. However the underlying concentration distribution is similar to that of the previous case. We see that the linear predictor is a good approximation when the concentration is near the reaction steady state value.

\begin{figure}[h]
\begin{center}
\includegraphics[width=1\linewidth]{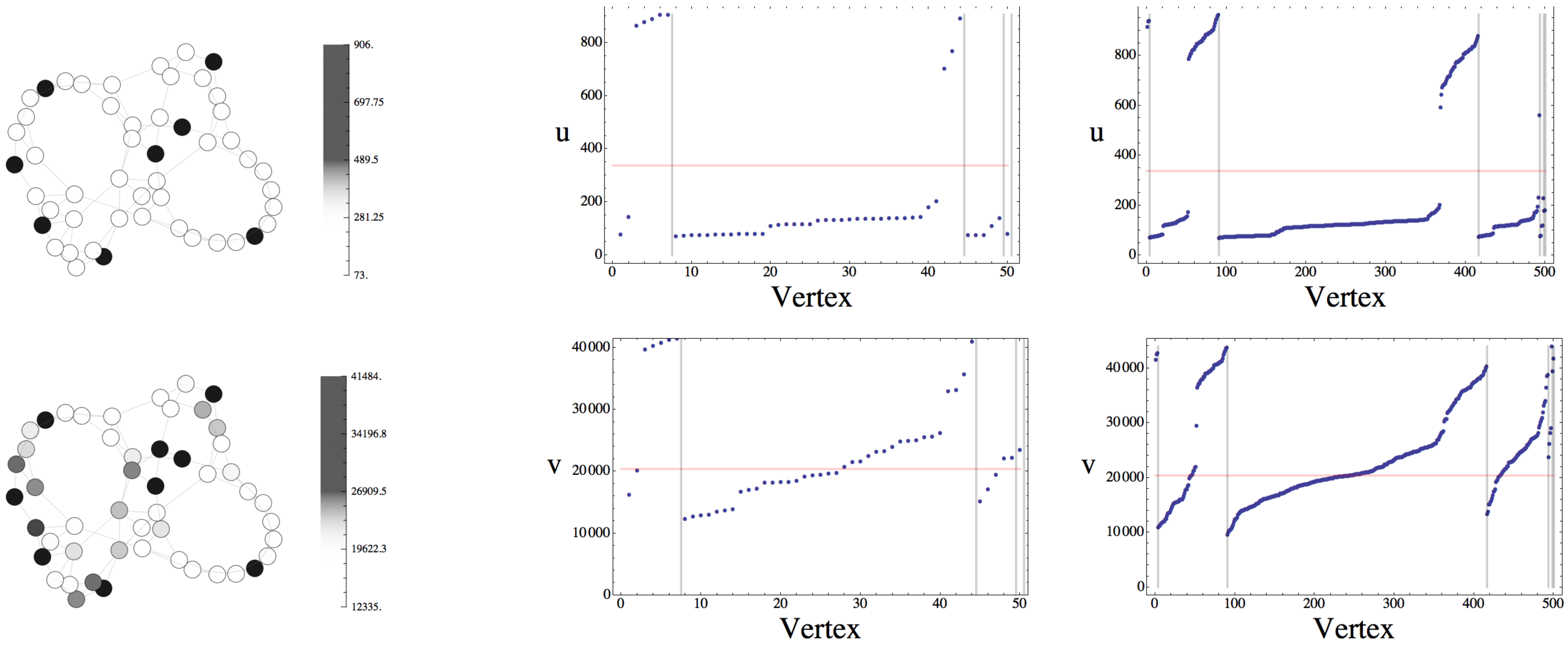}
\caption{ Case A Laplacian Gierer-Meinhardt reaction diffusion on a WS network with $k=3$ and $p=0.1$. The characteristic pattern of $u$ (\emph{top}) and $v$ (\emph{bottom}) on 50 (\emph{left}) vertex network and the corresponding distribution of concentrations on a network with 50 vertices (\emph{centre}) and 500 (\emph{right}) vertices. The horizontal line gives the reaction steady state.}
\label{fig:WS_500_GM_CA}
\end{center}
\end{figure}

\begin{figure}[h]
\begin{center}
\includegraphics[width=1\linewidth]{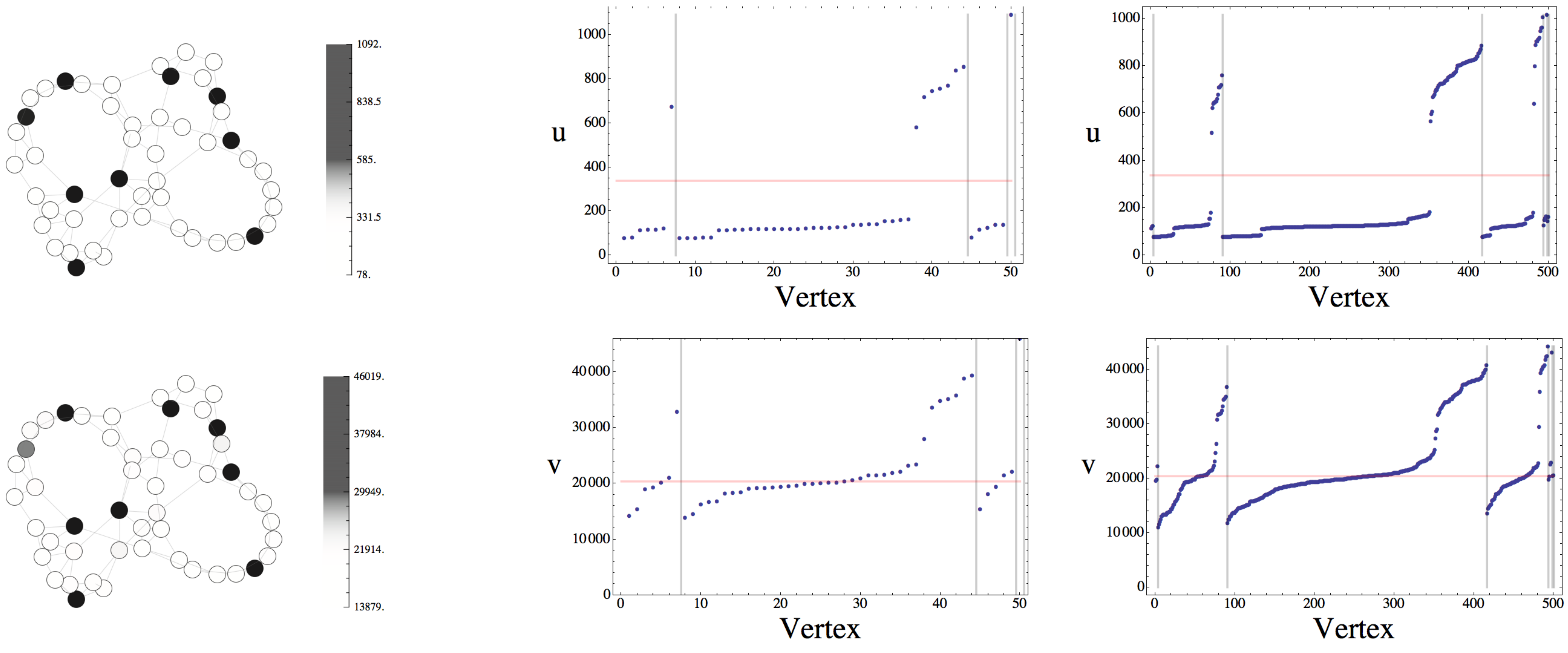}
\caption{ Case B Laplacian Gierer-Meinhardt reaction diffusion on a WS network with $k=3$ and $p=0.1$. The characteristic pattern of $u$ (\emph{top}) and $v$ (\emph{bottom}) on 50 (\emph{left}) vertex network and the corresponding distribution of concentrations on a network with 50 vertices (\emph{centre}) and 500 (\emph{right}) vertices. The horizontal line gives the reaction steady state.}
\label{fig:WS_500_GM_CB}
\end{center}
\end{figure}

We consider the Gierer-Meinhardt reaction kinetics with both Case A and Case B Laplacian operators as shown in Figs. \ref{fig:WS_500_GM_CA} and \ref{fig:WS_500_GM_CB} respectively. In both cases, the distribution of concentrations is bimodal for each vertex degree.

The example using the Case B Laplacian exhibits patterns very similar to those using the Case A. Similarly to the BA network, the highest concentrations in the second example are always increasing in contrast to the first example. Unlike the patterns on the BA network, the concentrations for both cases are at comparable quantities.

\begin{figure}[h]
\begin{center}
\includegraphics[width=1\linewidth]{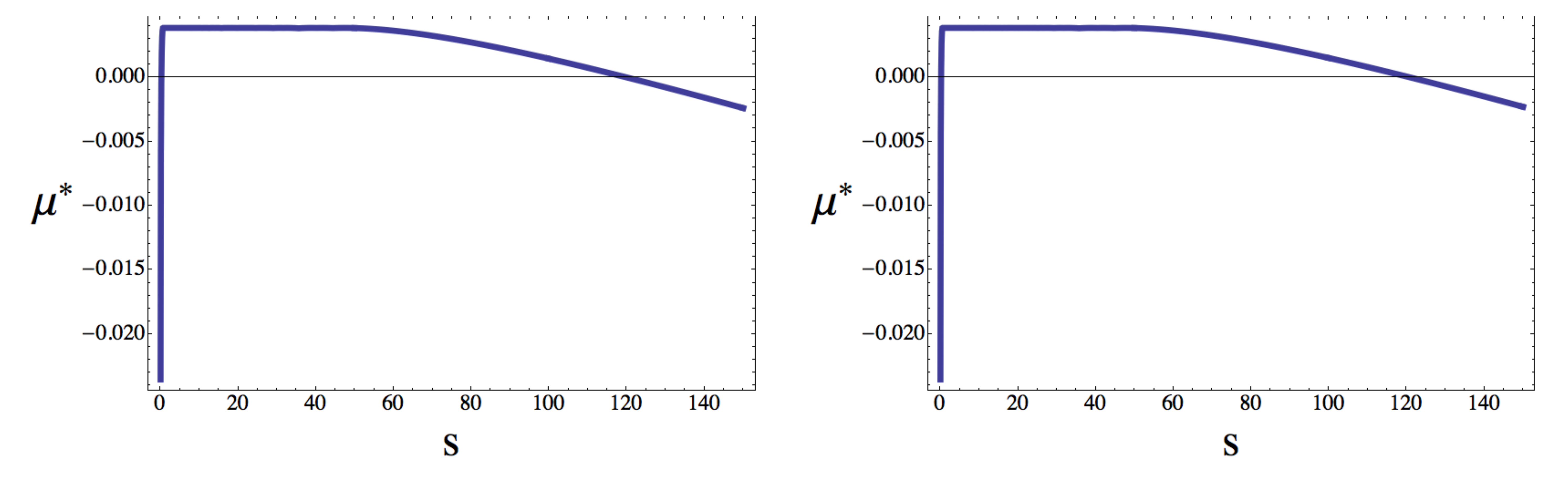}
\caption{Dispersion relation for Case A (\emph{left}) and Case B (\emph{right}) Laplacian Gierer-Meinhardt reaction diffusion on a 500 vertex WS network with $k=3$ and $p=0.1$.}
\label{fig:WS_500_GM_D}
\end{center}
\end{figure}

Dispersion relations for the Case A and Case B models are shown in Fig. \ref{fig:WS_500_GM_D}.
Unlike in the BA network, both Laplacian operators permit dispersion relations with almost indistinguishable profiles.

\section{Summary and Discussion}
\label{sec:conclusion}

In this paper we have derived the generalized master equation for reaction diffusion processes on networks based on the CTRW as the underlying stochastic process. We assumed that the vertices of the network can be occupied by many particles with the reactions occurring among the particles on the same vertex. The reactions were assumed to be governed by the same reaction kinetics on each vertex. The CTRW models particles jumping between vertices and reaction kinetics were incorporated into this as birth and death processes. We used the CTRW framework to derive a family of Laplacian operators that govern the diffusion on the network. These operators are dependent on the waiting time probability density at each vertex and the jumping probability density between vertices. In the case of non-exponential waiting time densities these operators convolve the reaction and transport processes. However in the case of exponential waiting times the operators are purely transport operators. The complication of power law waiting time densities with reactions on networks is an important problem for future work.

In general it is to be expected that the CTRW network reaction-diffusion model can lead to unequal concentrations of particles across the vertices. We investigated this pattern formation in the concentrations of particles. We considered the jumping probability to be equal across any edge from a given vertex and we considered the waiting time densities to be exponential, for two choices of the rate parameter: Case A, where the rate parameter was taken to be proportional to the vertex degree, and Case B, where the rate parameter was taken to be the same for all edges.

We identified three distinct pattern formation mechanisms in our CTRW network reaction-diffusion models. In the case of symmetric network Laplacian operators pattern formation followed a Turing mechanism whereby the steady state of the reaction dynamics, homogeneous across the vertices, became destabilized by the jumps resulting in a nonhomogeneous pattern across the vertices. The Turing mechanism can also result in pattern formation when the network Laplacian is nonsymmetric but other pattern formation mechanisms can occur in this case too. Second, the use of the nonsymmetric Laplacian may by itself, lead to a build up in concentrations on vertices according to their degree. Third, if the nonsymmetric Laplacian is coupled with reactions that have a finite nonzero steady state then the interplay between the two can result in a different pattern across the network. The different signatures of these patterns may prove useful in elucidating the underlying transport processes in networks when this is not known at the outset. 

The family of Laplacians we have derived should be used when the underlying process may be represented by a CTRW. This embodies a large class of diffusive processes on networks. By considering the underlying stochastic process, the \emph{ad hoc} choices of transport operators are constrained resulting in a physically consistent model of diffusion on networks. Our models of reaction-diffusion on networks are capable of reproducing a wide range of observed dynamics, as exemplified by our pattern formation examples.

\begin{acknowledgments}
This work was supported by the Australian Commonwealth Government (ARC No. DP1094680).
\end{acknowledgments}

\bibliography{network}

\end{document}